\documentclass[11pt]{article}
\usepackage[american]{babel}
\usepackage{amsmath}
\usepackage{latexsym}
\usepackage{amssymb,url,pstricks}
\usepackage{theorem}
\usepackage{diagrams}
\theoremstyle{change}
\newtheorem{Thm}{Theorem}[section]

\newtheorem{Lem}[Thm]{Lemma}
\newtheorem{Fact}[Thm]{Fact}
{\theorembodyfont{\rmfamily}

\newtheorem{Rem}[Thm]{Remark}
}
\renewcommand{\phi}{\varphi}
\renewcommand{\rho}{\varrho}

\newcommand{\proof}{\par\medskip\rm\emph{Proof. }}

\newcommand{\qed}{\ \hglue 0pt plus 1filll $\Box$}

\newcommand{\RR}{\mathbb{R}}

\newcommand{\SKIP}[1]{}

\DeclareMathAlphabet{\varcal}{U}{rsfs}{m}{it}
\begin{document}

\title{{\bf Strongly transitive actions on Euclidean buildings}}
\author{Linus Kramer\thanks{Supported by SFB 878}
{\ and }
Jeroen Schillewaert\thanks{Research supported by the Alexander Von Humboldt Foundation through a fellowship for experienced researchers}}
\date{}
\maketitle

\begin{abstract}
We prove a decomposition result for a group $G$ acting strongly transitively on the Tits boundary of a Euclidean building.
As an application we provide a local to global result for discrete Euclidean buildings, which generalizes results in the 
locally compact case by Caprace--Ciobotaru and Burger--Mozes. 
Let $X$ be a Euclidean building without cone factors. If a group $G$ of automorphisms of $X$ acts strongly transitively on
the spherical building at infinity $\partial X$, then the $G$-stabilizer of every affine apartment in $X$ contains all reflections along thick walls.
In particular $G$ acts strongly transitively on $X$ if $X$ is simplicial and thick.
\end{abstract}

\maketitle

\section*{Introduction}

We study strongly transitive actions on not necessarily discrete Euclidean buildings and prove a decomposition result.
As an application  we obtain the following local to global result for discrete Euclidean buildings. 

\medskip\noindent\textbf{Theorem A}
{\em Let $X$ be a thick simplicial Euclidean building and let $\partial X$ be the spherical building at infinity 
with respect to the complete apartment system of $X$. Let $G$ be a group of type-preserving 
automorphisms of $X$. Then the action of $G$ on $\partial X$ is strongly transitive if and only if
the action of $G$ on $X$ is strongly transitive, with respect to the complete apartment system of $X$.}

\medskip
A special case of Theorem A was proved by Caprace and Ciobotaru assuming in addition that $X$ is locally finite and that 
$G$ is a closed subgroup of $\mathrm{Aut}(X)$ \cite[Theorem 1.1]{CC}. 
For trees, this result was proved by Burger and Mozes \cite[Lemma 3.1.1.]{BM}, and also in \cite[Corollary~3.6]{CC}. 
Recently Ciobotaru and Rousseau proved an analogon of this special case in the more general context of hovels \cite{CG}.
Theorem A is a consequence of the following general result about Euclidean buildings.

\medskip\noindent\textbf{Theorem B}
{\em 
Let $X$ be a Euclidean building and let $\partial X$ be its spherical building at infinity with respect to the
complete apartment system. Assume that no irreducible factor
of $X$ is a Euclidean cone and that $\partial X$ is thick.
Suppose that a group $G$ of automorphisms of $X$
acts strongly transitively on $\partial X$. Let $A\subseteq X$ be an affine apartment in the complete apartment system and 
let  $c\in\partial A$ be a chamber.
Then the $G$-stabilizer $B$
of $c$ splits as \[B=TV,\] where $T\subseteq B$ is the point-wise stabilizer of $\partial A$ and
$V=\{g\in B\mid g\text{ fixes some point }p\in A\}\unlhd B$. 

If $M\subseteq A$ is a thick wall, then there exists an element $g\in G$ that
leaves $A$ invariant and acts as a reflection along $M$ on $A$.}


\section{The set-up.} 

Throughout, $X$ is a Euclidean building, as in Tits \cite{Tits}. However, we use the somewhat more
general definition where $X$ is allowed to be reducible.
Parreau's work \cite{Par} is an excellent reference.
Some additional results that we use can be found in \cite{KW}, in particular in Section 4.
We emphasize that we do not assume that $X$ is irreducible,
discrete, or metrically complete.

An \emph{automorphism} of $X$ is an isometry whose composition with the coordinate charts in the given atlas
results in the same collection of charts \cite[Definition 2.5]{Par}. 
The Tits boundary of the CAT$(0)$ space $X$ is then the metric realization of a simplicial
spherical building $\partial X$, whose type is given by the spherical Coxeter group of the affine Weyl group
of $X$. 
The spherical apartments in $\partial X$ are in one-to-one correspondence with the
affine apartments in the complete apartment system \cite[Corollaire 2.19]{Par}. 
We note that the Tits boundary is usually defined only for complete CAT$(0)$ spaces. However,
its definition in terms of equivalence classes of geodesic rays works well for euclidean buildings
even if they are not complete \cite{Par}.

If the spherical building $\partial X$ is thick, then every isometry of $X$ induces a 
(not necessarily type-preserving) automorphism of $\partial X$. It follows then that
the isometry group $\mathrm{Iso}(X)$ has a subgroup $G$ of finite index which
consists of automorphisms of $X$ with respect to the complete apartment system.
So the difference between isometries of $X$ and automorphisms of $X$ for the complete apartment system
is not very big in this situation.

We briefly recall what strong transitivity means \cite[6.1.1]{AB}. A group $G$ of type-preserving automorphisms
acts \emph{strongly transitively} on a (combinatorial) building $\Delta$ if $G$ acts transitively on the apartments
of $\Delta$ (for a given apartment system),
and if the set-wise stabilizer of an apartment acts transitively on the chambers in the apartment.
If $N\subseteq G$ is the set-wise $G$-stabilizer of an apartment $\Sigma$ and if
$c\in\Sigma$ is a chamber, with $G$-stabilizer $B$, then $(G,B,N)$ constitutes a \emph{Tits system}
(or BN-pair) for $G$ \cite[6.2.6]{AB}. For spherical buildings, strong transitivity may also be stated as follows:
$G$ acts transitively on the chambers of $\Delta$, and the $G$-stabilizer $B$ of a chamber
$c$ acts transitively on all chambers opposite $c$ \cite[Proposition 6.15]{AB}.

The \emph{affine Weyl group} of $X$ is the semidirect product of a spherical Coxeter group $W_0$ and
a $W_0$-invariant translation group $T\subseteq\RR^m$, with $T\otimes\mathbb R=\mathbb R^m$.
\emph{Walls}, \emph{half-apartments} and \emph{Weyl chambers} in $X$ are defined as in \cite{Par}
(a Weyl chamber is what Tits \cite{Tits} calls a \emph{sector}).
A wall is called \emph{thick} if it can be written as the intersection of three affine apartments.
The visualizations below of a Weyl chamber, an affine half-apartment and a wall may be helpful,
see \cite[Section 4]{KW}.
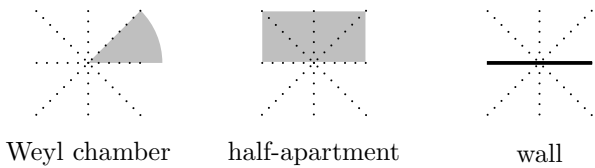
\begin{center}
\begin{pspicture}(9,2)
\psline[linewidth=0.5mm](6.3,1)(7.7,1)
\psframe[linestyle=none,fillstyle=solid,fillcolor=lightgray](3.3,1)(4.7,1.7)
\pswedge[linestyle=none,fillstyle=solid,fillcolor=lightgray](1,1){1}{0}{45}
\psline[linestyle=dotted](.3,.3)(1.7,1.7)
\psline[linestyle=dotted](.3,1)(1.7,1)
\psline[linestyle=dotted](1,.3)(1,1.7)
\psline[linestyle=dotted](.3,1.7)(1.7,.3)
\psline[linestyle=dotted](3.3,.3)(4.7,1.7)
\psline[linestyle=dotted](3.3,1)(4.7,1)
\psline[linestyle=dotted](4,.3)(4,1.7)
\psline[linestyle=dotted](3.3,1.7)(4.7,.3)
%
\psline[linestyle=dotted](6.3,.3)(7.7,1.7)
\psline[linestyle=dotted](6.3,1)(7.7,1)
\psline[linestyle=dotted](7,.3)(7,1.7)
\psline[linestyle=dotted](6.3,1.7)(7.7,.3)
\rput(1,-.2){\small Weyl chamber}
\rput(4,-.2){\small half-apartment}
\rput(7,-.2){\small wall}
\end{pspicture}\\\
\end{center}
For an affine apartment $A$ or an affine half-apartment $H$ we denote
by $\partial A\subseteq\partial X$ and $\partial H\subseteq\partial X$ the set of all Weyl simplices
at infinity which have representatives in $A$ and $H$, respectively.

\begin{Fact}\label{Facts}
The following facts can be found in \cite[Corollaire 1.6 and Proposition 1.7]{Par} and \cite[Lemma 5.14]{KW}.

(i) The affine apartments in $X$ (with respect to the complete apartment system) are in one-to-one correspondence with the apartments
in $\partial X$ via
$A\longmapsto\partial A$.

(ii) If $c\in\partial X$ is a chamber and if $A_1,\ldots,A_\ell\subseteq X$ are affine apartments
with $c\in\partial A_1\cap\cdots\cap \partial A_\ell$, then $A_1\cap\cdots \cap A_\ell$ is a closed convex set
which contains a Weyl chamber $S$ with $\partial S=c$.
\end{Fact}

\section{The proofs.}

Our convention for commutators is that $[a,b]=ab(ba)^{-1}$, and we let maps act from the left.
The proof of Theorem B relies now on the following lemmata. 
\begin{Lem}\label{TelescopeApartment}
Let $X$ be a Euclidean building and let $H\subseteq X$ be an affine half-apartment.
Suppose that $g$ is an automorphism of $X$ that fixes $\partial H$ point-wise. If $H\subsetneq g(H)$, then
$E=\bigcup_{n\geq 0}g^n(H)$ is an affine apartment in the complete apartment system.
Moreover, $E$ is invariant under the cyclic group $\langle g\rangle$ generated by $g$.
The element $g$ acts as a translation on $E$.

\proof
Let $M\subseteq H$ denote the wall that bounds $H$. Then $g(M)$ is parallel to $M$,
and we denote the Hausdorff distance between $M$ and $g(M)$ by $\lambda >0$.
We fix an isometric embedding $\phi_0:H\to\mathbb R^m$, where $m=\dim(H)$, such that
$H$ maps onto $(-\infty,0]\times\mathbb R^{m-1}$. Then $\phi_0$ extends uniquely to an
isometric embedding $\phi_1:g(H)\to\mathbb R^m$, with image $(-\infty,\lambda]\times\mathbb R^m$.
Continuing in this fashion, we obtain isometries 
$\phi_n:g^n(H)\to (-\infty,n\lambda]\times\mathbb R^{m-1}$ which fit together to a
map $\phi_\infty:E\to\mathbb R^m$. From the construction it is clear that $\phi_\infty$
is also an isometric embedding, and since $\mathbb R=\bigcup_{n\geq 0}(\infty,n\lambda]$
the map $\phi_\infty$ is onto.
By \cite[Proposition 2.25]{Par}, $E$ is an affine apartment in the complete apartment
system. From the construction of $E$ it is clear that $g(E)=E$,
hence $E$ is invariant under $\langle g\rangle$.
Finally we note that $g$ fixes $\partial H\subseteq\partial E$ point-wise.
It follows that $g$ fixes $\partial E$ point-wise, and hence that $g$ acts as a translation
on $E$.
\qed
\end{Lem}
\begin{Lem}\label{TranslationLemma}
Suppose that $X$ is a Euclidean building and that $g$ is an automorphism of $X$ that fixes
a chamber at infinity $c\in\partial X$. Suppose that $A\subseteq X$ is an affine apartment
with $c\in\partial A$. Then there exists a unique translation $t$ in the isometry group of $A$
such that $g^{-1}(p)=t(p)$ holds for all $p\in A\cap g(A)$. 

\proof
The convex set $Q=g(A)\cap A$ contains a
Weyl chamber $S$ by \ref{Facts}(ii). Now $g^{-1}$ restricts to an isometric embedding $t:Q\to A$.
Since $Q\subseteq A$ is convex and top-dimensional, there is a unique isometric extension
of $t$ to an isometry $t:A\to A$. This isometry is given by an element of the affine
Weyl group, because $g$ is an automorphism and because the transition maps between
apartments are given by elements in the affine Weyl group. 
Hence it is a product of a translation and an element in the spherical
Coxeter group $W_0$. Since $t$ fixes the chamber $c$, the map $t$ is a translation.
\qed
\end{Lem}
The proof of the next Lemma follows closely the reasoning in \cite[Proposition 2.3]{KW}.
\begin{Lem}\label{AlmostThere}
Let $X$ be a Euclidean building and let $A\subseteq X$ be an affine apartment.
Let $M\subseteq A$ be a thick wall. Let $H_1,H_2\subseteq A$ denote the two
affine half-apartments bounded by $M$. Let $H_3\subseteq X$ be a third affine
half-apartment, with $M=H_1\cap H_3=H_2\cap H_3$. Let $K$ denote a
group of automorphisms of $X$ fixing $\partial H_1$ point-wise.
If $K$ acts transitively on the set of all spherical apartments in $\partial X$
containing $\partial H_1$, then there exists an element $k\in K$ that
fixes $H_1$ set-wise and maps $H_2$ to $H_3$. In particular, $k$ fixes the wall
$M$ set-wise.

If in addition $K$ contains an element $t$ that acts as a translation on $A$
which does not leave $M$ invariant, then the element $k$ above 
can be chosen in such a way that 
$k$ fixes $M$ point-wise.

\proof
Let $g\in K$ be an element that maps $\partial H_2$ to $\partial H_3$.
Such an element exists by assumption. There are three cases: $g(H_1)=H_1$,
or $g(H_1)\supsetneq H_1$, or $g(H_1)\subsetneq H_1$.
In the first case $k=g$ fixes $M$ set-wise and we are done.
In the second case we use \ref{TelescopeApartment}. There exists a $\langle g\rangle$-invariant
affine apartment $E$ with $H_1\subseteq E$, and $g$ acts as a translation on $E$.
Let $h\in K$ be an element that maps $E$ to $A$. Such an element exists by our assumptions.
Then $[g,h]$ maps $A$ to $H_1\cup H_3$. Applying \ref{TranslationLemma} to the four automorphisms
$g,h,g^{-1}$ and $h^{-1}$ we find four
Weyl chambers $S_1,S_2,S_3,S_4\subseteq A$ with 
$\partial S_1=\partial S_2=\partial S_3=\partial S_4\subseteq\partial H_1$
on which $g,g^{-1},h,h^{-1}$ act as isometric embeddings $S_i\rTo A$ which are restrictions of translations $A\rTo A$.
Passing to subsectors, we may assume that we have both
$S_4\subseteq S_3\subseteq S_2\subseteq S_1$ and $S_4\rTo^{h^{-1}} S_3\rTo^{g^{-1}}S_2\rTo^{h} S_1\rTo^g A$.
Since $g$ and $h$ are restrictions of translations
and since translations commute,
$[g,h]$ is the identity on $S_4$. Therefore $[g,h]$ fixes the wall $M$ point-wise.
The remaining case is similar, with $g$ replaced by $g^{-1}$.

If $K$ contains an element $t$ that acts as a translation on $A$ and does not leave $M$
invariant and if $g(H_1)=H_1$, then $gt$ still maps $\partial H_2$ to $\partial H_3$,
but $gt(H_1)\neq H_1$. The construction above now gives us an element $k$ that fixes
$M$ point-wise.
\qed
\end{Lem}

\begin{Lem}\label{ShiftLemma}
Let $X$ be a Euclidean building and suppose that a group $G$ of automorphisms of $X$ acts
strongly transitively on $\partial X$. Let $A\subseteq X$ be an affine apartment and
suppose that $A$ contains two distinct parallel thick walls $M_1,M_2$.
Then there exists an element $t\in G$ that fixes $\partial A$ point-wise and that
maps $M_1$ onto its image under the reflection along $M_2$.

\proof
For $i=1,2$, let $H_i\subseteq X$ be an affine half-apartment bounded by $M_i$, with 
$H_i\cap A=M_i$. Denote the affine half-apartments in $A$
bounded by $M_1$ by $H_+$ and $H_-$.
Applying \ref{AlmostThere} three times to the three affine apartments $A$, $H_1\cup H_+$ and $H_1 \cup H_-$, 
we see that there exists an element $g_1\in G$ that maps $A$ onto itself,
fixes $M_1$ set-wise, and interchanges $H_+$ and $H_-$. Similarly, we find an element $g_2$ fixing $M_2$ set-wise
and interchanging the two affine half-apartments in $A$ bounded by $M_2$.
We note that $g_1$ and $g_2$ act in the same way on $\partial A$, namely as reflections
along the wall $\partial M_1=\partial M_2$.
Then $g_2g_1$ is a translation and maps $M_1$ to the reflection of $M_1$ along $M_2$.
\qed
\end{Lem}
The next result about spherical buildings is folklore.
\begin{Lem}\label{GenerationLemma}
Suppose that a group $G$ acts strongly transitively on a thick spherical building $\Delta$. Let $\Sigma\subseteq\Delta$ be an apartment
containing a chamber $c$, and let $B$ denote the $G$-stabilizer of $c$. Let $T$ denote the point-wise $G$-stabilizer of $\Sigma$.
For each pair of opposite panels $\tau=(a,\bar a)$ in $\Sigma$, let $B_\tau$ denote the point-wise $B$-stabilizer of $\tau$.
Then $B$ is generated by $T$ and the groups $B_\tau$, where $\tau$ runs over all pairs of opposite panels in $\Sigma$.

\proof
Let $K\subseteq B$ denote the group generated by $T$ and the $B_\tau$.
Let $\bar c\in\Sigma$ denote the unique chamber opposite $c$ in $\Sigma$
and let $c=c_0,\ldots,c_m=\bar c$ be a minimal gallery in $\Sigma$.
Let $b\in B$ and let $i\geq 0$ be the largest index such that $b(c_i)\in\Sigma$. We prove by induction on $m-i$ that
$b\in K$. For $m=i$ we have that $b\in T\subseteq K$. For $m>i$,
the intersection $c_i \cap c_{i+1}$ is a panel in $\Sigma$, which has a unique opposite in $\Sigma$. So $c_i$ and $c_{i+1}$
determine a pair $\tau$ of opposite panels in $\Sigma$. By strong transitivity, there exists an
element $k\in B_\tau\subseteq K$ with $k(c_{i+1})=b(c_{i+1})$, and $k(c_j)=c_j=b(c_j)$ for all $j\leq i$. 
Then $k^{-1}b\in K$ by the induction hypothesis, and thus $b\in K$.
\qed
\end{Lem}
\emph{Proof of Theorem B.}
If $b\in B$ fixes a point $p\in A$, then $b$ fixes a Weyl chamber $S\subseteq A$ with tip $p$
and with $\partial S=c$. Hence $V$ consists of all elements in $B$ which fix a
Weyl chamber $S\subseteq A$ with $\partial S=c$. It follows from \ref{Facts}(ii) that $V$
is a subgroup of $B$, and clearly $V\unlhd B$.

Suppose that $b\in B$ is an element such that $b(A)\cap A=H$ is an affine half-apartment.
Let $M$ denote the thick wall bounding $H$. Since none of the irreducible factors of $X$ is a
Euclidean cone, there exist thick walls in $A$ different from $M$ which are parallel to
$M$, see \cite[Proposition 4.21]{KW}. Using \ref{ShiftLemma} and then \ref{AlmostThere},
there exists $v\in V$ with $b(A)=v(A)$, whence $b\in VT$. By \ref{GenerationLemma},
$B$ is generated by $T$ and such elements $b$. It follows that $B=VT$.

Given a thick wall $M\subseteq A$ bounding affine half-apartments $H_1,H_2\subseteq A$ and
a third affine half-apartment $H_3$ bounded by $M$
with $H_3 \cap A=M$, we can therefore find elements $g_1,g_2,g_3\in G$
fixing $M$ point-wise, with $A\xrightarrow{g_1} H_1\cup H_3\xrightarrow{g_2} H_3\cup H_2\xrightarrow{g_3} A$.
Then $g_3g_2g_1\in N$ fixes $M$ point-wise and interchanges $H_1$ and $H_2$.
\qed

\medskip\noindent
\emph{Proof of Theorem A.}
Suppose that $G$ acts strongly transitively on $\partial X$. Since the affine apartments in the complete apartment system are in
one-to-one correspondence with the apartments in $\partial X$, the group $G$ acts transitively on the complete apartment system.
Let $A\subseteq X$ be an affine apartment and let $C\subseteq A$ be a chamber. Each codimension $1$ face of $C$ determines
a unique wall $M$ in $A$, which is thick because $X$ is thick. By Theorem B, the group $G$ contains an element $g$
which leaves $A$ invariant and acts on $A$ as a reflection along $M$. These reflections for the faces of $C$ generate a 
group  which acts as the Coxeter group of $X$ on $A$; in particularly, it acts transitively on the chambers in $A$.
Thus $G$ acts strongly transitively on $X$.

Conversely, suppose that $G$ acts strongly transitively on $X$. Then $G$ acts transitively on the apartments of
$\partial X$. Let $A\subseteq X$ be an affine apartment and let $v\in A$ be a special vertex. Each Weyl chamber
$S\subseteq A$ with tip $v$ contains a unique chamber $C$ with $v\in C$. Since the group $N\subseteq G$
consisting of all elements $g\in G$ that leave $A$ invariant acts transitively on these chambers in $A$ containing $C$,
it acts transitively on the Weyl chambers in $A$ with tip $v$. It follows that $N$ acts transitively on the chambers
in $\partial A$.
\qed
\begin{Rem}
The splitting $B=VT$ with $[B,B]\subseteq V$ in Theorem B looks deceivingly 
similar to a split Tits system for $\partial X$.
However, this is misleading. The group $V$ need not be nilpotent or solvable, as the case of the
Bruhat--Tits building for $\mathrm{PGL}_{m+1}(D)$ shows, where $D$ is a central division ring over $\mathbb Q_p$
and $m\geq 1$. In this case $V$ contains the norm $1$ group $[D^*,D^*]=\mathrm{SL}_1(D)$, which is not solvable~\cite{H}. 

We observe also that the conclusion of the result holds, rather trivially, for Euclidean cones over
spherical buildings (vectorial buildings). However, we do not know if our result holds 
if $X$ is a nontrivial product of Euclidean cones and 'honest' Euclidean buildings.

Lastly, note that we always worked with the complete apartment system of the Euclidean building,
and this was crucial in the proof of \ref{TelescopeApartment}.
We have no results for strongly transitive actions with respect to smaller apartment systems.
\end{Rem}

{{\bigskip
\raggedright
Linus Kramer, Jeroen Schillewaert\\
Mathematisches Institut, 
Universit\"at M\"unster,
Einsteinstr. 62,
48149 M\"unster,
Germany\\
{\makeatletter
e-mail: {\tt linus.kramer{@}wwu.de, jschillewaert@{g}mail.com}}}
 
\end{document}